\documentclass{article}
\usepackage{amsmath}

\newenvironment{proof}[1][Proof]{\noindent\textbf{#1.} }{\ \rule{0.5em}{0.5em}}
\input{tcilatex}

\begin{document}

\section{Introduction}

\qquad Let $R$\ be a convex quadrilateral in the $xy$ plane. An ellipse
which passes through the vertices of $R$\ is called a circumscribed ellipse
or ellipse of circumscription. In the book ([D]), Dorrie presents Steiner's
nice characterization of the ellipse of circumscription which has minimal
eccentricity, which he calls the most nearly circular ellipse. A pair of 
\textbf{conjugate diameters} are two diameters of an ellipse such that each
bisects all chords drawn parallel to the other. Let $\theta _{1}$ and $%
\theta _{2}$ be the angles which a pair of conjugate diameters make with the
positive $x$ axis. Then $\tan \theta _{1}$ and $\tan \theta _{2}$ are called
a pair of \textbf{conjugate directions}. First, Steiner proves that there is
only one pair of conjugate directions, $M_{1}$ and $M_{2}$, that belong to
all ellipses of circumscription. Then he proves in essence that \textbf{if}
there is an ellipse, $E$, whose \textbf{equal} conjugate diameters possess
the directional constants $M_{1}$ and $M_{2}$, then $E$ must be an ellipse
of circumscription which has minimal eccentricity. There are several gaps
and missing pieces in Steiner's result. Steiner does \textbf{not }show that
there \textbf{exists} an ellipse of circumscription, $E$, whose equal
conjugate diameters possess the directional constants $M_{1}$ and $M_{2}$,
or that such an ellipse is \textbf{unique}. \ He also does \textbf{not }prove%
\textbf{\ }in general the \textbf{uniqueness} of an ellipse of
circumscription which has minimal eccentricity . It is possible that there
could exist a circumscribed ellipse of minimal eccentricity that might 
\textbf{not} have \textbf{equal } conjugate diameters which possess the
directional constants $M_{1}$ and $M_{2}$.

In Propositions 1 and 2 we fill in these gaps in Steiner's proof. We
prove(Proposition 1), without using the directional constants $M_{1}$ and $%
M_{2}$, that there is a unique ellipse of minimal eccentricity which passes
through the vertices of $R$. Then we show(Proposition 2) that there exists
an ellipse which passes through the vertices of $R$\ and whose \textit{equal}
conjugate diameters possess the directional constants $M_{1}$ and $M_{2}$.
In addition, our methods enable us to prove(Theorem 2) that there is a
unique ellipse of \textbf{minimal area} which passes through the vertices of 
$R$.

In [H] the author proved numerous results about ellipses \textbf{inscribed }%
in convex quadrilaterals. We filled in similar gaps in a problem often
referred to in the literature as Newton's problem, which was to determine
the locus of centers of ellipses inscribed in $R$. In particular, in [H] the
author proved that there exists a unique ellipse of minimal eccentricity, $%
E_{I}$, inscribed in $R$. This leads to the last section of this paper,
where we discuss a special class of convex quadrilaterals which we call
bielliptic and which generalize the bicentric quadrilaterals. A convex
quadrilateral, $R$, is called bicentric if there exists a circle inscribed
in $R$\ and a circle circumscribed about $R$. Let $E_{O}$ denote the unique
ellipse of minimal eccentricity circumscribed about $R$. $R$\ is called 
\textbf{bielliptic} if $E_{I}$ and $E_{O}$ have the \textbf{same}
eccentricity. We prove(Theorem 4), that there exists a bielliptic convex
quadrilateral which is not bicentric. We also prove(Theorem 5), that there
exists a bielliptic trapezoid which is not bicentric.

Finally we prove the perhaps not so obvious result(Theorem 3), that if $E_{1}
$ and $E_{2}$ are each ellipses, with $E_{1}$ inscribed in $R$\ and $E_{2}$
circumscribed about $R$, then $E_{1}$ and $E_{2}$ cannot have the same
center.

\section{Minimal Eccentricity}

We state the following lemma without proof. The details are well--known and
can be found in

numerous places.

\textbf{Lemma 1:} Consider the ellipse, $E_{0}$, with equation $%
Ax^{2}+By^{2}+2Cxy+Dx+Ey+F=0,A,B>0,AB-C^{2}>0,AE^{2}+BD^{2}+4FC^{2}-2CDE-4ABF\neq 0 
$; Let $a$ and $b$ denote the lengths of the semi--major and semi--minor
axes, respectively, of $E_{0}$. Let $\phi $ denote the acute rotation angle
of the axes of $E_{0}$ going counterclockwise from the positive $x$ axis and
let $\left( x_{0},y_{0}\right) $ denote the center of $E_{0}$. Then

\begin{equation}
a^{2}=\dfrac{AE^{2}+BD^{2}+4FC^{2}-2CDE-4ABF}{2(AB-C^{2})\left( (A+B)-\sqrt{%
(B-A)^{2}+4C^{2}}\right) },  \tag{a}
\end{equation}

\begin{equation}
b^{2}=\dfrac{AE^{2}+BD^{2}+4FC^{2}-2CDE-4ABF}{2(AB-C^{2})\left( (A+B)+\sqrt{%
(B-A)^{2}+4C^{2}}\right) },  \tag{b}
\end{equation}

\begin{equation}
\phi =\dfrac{1}{2}\cot ^{-1}\left( \dfrac{A-B}{2C}\right) ,C\neq 0\text{ and 
}\phi =0\text{ if }C=0,  \tag{phi}
\end{equation}

and 
\begin{equation}
x_{0}=\allowbreak -\dfrac{1}{2}\dfrac{BD-CE}{AB-C^{2}},y_{0}=\allowbreak 
\dfrac{1}{2}\dfrac{CD-AE}{AB-C^{2}}  \tag{center}
\end{equation}

Throughout this section, we let $R$\ be a given convex quadrilateral in the $%
xy$ plane, and we assume throughout that $R$\ is not a trapezoid. We use the
notation and terminology used by Steiner (see [D]), and we also assume that $%
R$\ has the shape given in [D]. Other shapes for a convex quadrilateral are
possible, of course, but we do not consider those cases in the proofs below,
the details being similar. Let $OPRQ$ denote the vertices of $R$, in
counterclockwise order. By using an isometry of the plane, we can assume
that $O=(0,0)$ and that $P,R,$ and $Q$ are in the first quadrant. Let $H=%
\overleftrightarrow{QR}\cap \overleftrightarrow{OP},K=\overleftrightarrow{PR}%
\cap \overleftrightarrow{OQ},q=\left\vert \overline{OQ}\right\vert
,h=\left\vert \overline{OH}\right\vert ,$ and $k=\left\vert \overline{OK}%
\right\vert $. Use the oblique coordinate system with $\overrightarrow{OP}$
as the $x$ axis and $\overrightarrow{OQ}$ as the $y$ axis, with those sides
given by $y=0$ and $x=0$. The sides $\overleftrightarrow{PR}$ and $%
\overleftrightarrow{QR}$ are given by $z=0$ and $w=0,$ respectively, where $%
z=kx+py-kp$ and $w=qx+hy-hq$. It follows that 
\begin{equation}
0<p<h,0<q<k.  \tag{1}
\end{equation}

Any ellipse passing through the vertices of $R$\ has equation $\lambda
xz+\mu yw=0$, where $\lambda $ and $\mu $ are \textit{nonzero} real numbers.
Letting $v=\dfrac{\lambda }{\mu }$, the equation becomes $vxz+yw=0,$ or

\begin{equation}
kvx^{2}+hy^{2}+(vp+q)xy-vkpx-hqy=0.  \tag{ell}
\end{equation}

Using $A=kv$, $B=h$, $C=\dfrac{1}{2}(vp+q)$, $D=-vkp$, $E=-hq$, and $F=0,$
(ell) is the equation of a nontrivial ellipse(that is, not just a single
point) if and only if $AB-C^{2}=khv-\dfrac{1}{4}(vp+q)^{2}>0$ and $%
AE^{2}+BD^{2}+4FC^{2}-2CDE-4ABF=khv\left( vp^{2}(k-q)+q^{2}(h-p)\right) \neq
0$. Note that $khv>0$ implies that $v>0$, so that $khv\left(
vp^{2}(k-q)+q^{2}(h-p)\right) >0$. We write the first condition, $AB-C^{2}>0$%
, as

\begin{equation*}
g(v)>0\text{, }g(v)=-p^{2}v^{2}+\left( 4kh-2pq\right) v-q^{2}.
\end{equation*}

Note that $g(v)=0\iff v=\dfrac{1}{p^{2}}\left( 2kh-pq\pm 2\sqrt{kh\left(
kh-pq\right) }\right) $. Hence $g(v)>0\iff v\in I$, where

$I=\left( \dfrac{1}{p^{2}}\left( 2kh-pq-2\sqrt{kh\left( kh-pq\right) }%
\right) ,\dfrac{1}{p^{2}}\left( 2kh-pq+2\sqrt{kh\left( kh-pq\right) }\right)
\right) $. Note that $\left( 2kh-pq\right) ^{2}-4\left( kh\left(
kh-pq\right) \right) =\allowbreak q^{2}p^{2}>0\Rightarrow 2kh-pq>2\sqrt{%
kh\left( kh-pq\right) }$ since $kh-pq>0$. Hence $I\subset \left( 0,\infty
\right) $. Our first main result is the following.

\textbf{Proposition 1:} There is a unique ellipse, $E_{O}$, of minimal
eccentricity which passes through the vertices of $R$.

\textbf{Proof: }By Lemma 1,

\begin{equation}
a^{2}=\dfrac{2khv\left( vp^{2}(k-q)+q^{2}(h-p)\right) }{\left( 4khv-\left(
vp+q\right) ^{2}\right) \left( (kv+h)-\sqrt{(kv-h)^{2}+\left( vp+q\right)
^{2}}\right) }  \tag{asq}
\end{equation}%
and 
\begin{equation}
b^{2}=\dfrac{2khv\left( vp^{2}(k-q)+q^{2}(h-p)\right) }{\left( 4khv-\left(
vp+q\right) ^{2}\right) \left( (kv+h)+\sqrt{(kv-h)^{2}+\left( vp+q\right)
^{2}}\right) },  \tag{bsq}
\end{equation}

which implies that $\dfrac{b^{2}}{a^{2}}=\dfrac{kv+h-\sqrt{(kv-h)^{2}+\left(
vp+q\right) ^{2}}}{kv+h+\sqrt{(kv-h)^{2}+\left( vp+q\right) ^{2}}}$. Some
simplification yields

\begin{equation}
\dfrac{b^{2}}{a^{2}}=f(v)=\dfrac{g(v)}{\left( (kv+h)+\sqrt{(kv-h)^{2}+\left(
vp+q\right) ^{2}}\right) ^{2}}.  \tag{fv}
\end{equation}

\ We shall now minimize the eccentricity by maximizing\textbf{\ }$\dfrac{%
b^{2}}{a^{2}}$. Now $f\,^{\prime }(v)=-2\dfrac{\left( 2hk-pq\right) \left(
vk-h\right) +p^{2}hv-q^{2}k}{\sqrt{(kv-h)^{2}+\left( vp+q\right) ^{2}}\left(
kv+h+\sqrt{(kv-h)^{2}+\left( vp+q\right) ^{2}}\right) ^{2}}$. Thus $%
f\,^{\prime }(v)=0\iff \left( 2hk-pq\right) \left( vk-h\right)
+p^{2}hv-q^{2}k=0\iff v=v_{0}$, where 
\begin{equation}
v_{0}=\dfrac{q^{2}k+2kh^{2}-hpq}{2k^{2}h-kpq+hp^{2}}.  \tag{v0}
\end{equation}%
Some simplification yields $(kv_{0}-h)^{2}+\left( v_{0}p+q\right) ^{2}=%
\dfrac{\left( ph+qk\right) ^{2}W}{\left( 2k^{2}h-kpq+hp^{2}\right) ^{2}}$,
which implies that 
\begin{equation}
g\left( v_{0}\right) =\dfrac{4kh\left( kh-pq\right) W}{\left(
2k^{2}h-kpq+hp^{2}\right) ^{2}},  \tag{gv0}
\end{equation}%
where

\begin{equation}
W=4k^{2}h^{2}+\left( hp-qk\right) ^{2}.  \tag{W}
\end{equation}%
Thus $g\left( v_{0}\right) >0$ by (gv0), which implies that $v_{0}\in I$.
Note that $kv+h+\sqrt{(kv-h)^{2}+\left( vp+q\right) ^{2}}>0$ for all $v>0$,
and $g(v)>0,v\in I$. Thus $f$ is differentiable on $I$ and $f$ has a unique
real critical point in $I$. Since $g$ vanishes at the endpoints of $I,$ $f$
also vanishes at the endpoints of $I$ since the denominator in (fv) is
positive at the endpoints of $I$. Since $f(v)>0$ on $I$, $f(v_{0})$ must
give the unique maximum of $f$ on $I$. 
\endproof%

In [D], Steiner shows that the unique pair of conjugate directions that
belong to all ellipses which pass through the vertices of $R$\ is given by

\begin{equation}
M_{1}=-\dfrac{k}{p}+\dfrac{r}{hp},M_{2}=-\dfrac{k}{p}-\dfrac{r}{hp},\text{
where }r=\sqrt{hk}\sqrt{hk-pq}.  \tag{M1M2}
\end{equation}

\textbf{Proposition 2:} There exists an ellipse which passes through the
vertices of $R$\ and whose \textit{equal} conjugate diameters possess the
directional constants $M_{1}$ and $M_{2}$.

\textbf{Proof: }Let $E_{O}$ denote the the unique ellipse from Proposition 1
of minimal eccentricity which passes through the vertices of $R$, and let $L$
and $L^{\prime }$ denote a pair of \textit{equal} conjugate diameters of $%
E_{O}$ with directional constants $M$ and $M^{\prime }$. Let $\phi $ denote
the acute angle of counterclockwise rotation of the axes of $E_{O}$ and let $%
a$ and $b$ denote the lengths of the semi--major and semi--minor axes,
respectively, of $E_{0}$. It is known(see, for example [S]) that $L$ and $%
L^{\prime }$ make equal acute angles, on opposite sides, with the major axis
of $E_{O}$. Let $\theta $ denote the acute angle going counterclockwise from
the major axis of $E_{O}$ to one of the equal conjugate diameters, which
implies that $\tan \theta =\dfrac{b}{a}$. By Lemma 1, with $A=kv$, $B=h$, $C=%
\dfrac{1}{2}(vp+q)$, $D=-vkp$, $E=-hq$, and $F=0,\cot \phi =\dfrac{kv-h}{vp+q%
}$. As one would expect from the results in [D], \textit{if} there is some
ellipse whose equal conjugate diameters possess the directional constants $%
M_{1}$ and $M_{2}$, then that ellipse minimizes the eccentricity among all
ellipses of circumscription. By the proof of Proposition 1, the point $v_{0}$
given in (v0) yields the ellipse which minimizes the eccentricity. Thus, to
prove Proposition 2, we let $v=v_{0}$. Then $\cot \phi =\dfrac{kv_{0}-h}{%
v_{0}p+q}=\dfrac{kq-hp}{2kh}\Rightarrow \phi =\dfrac{1}{2}\cot ^{-1}\left( 
\dfrac{kq-hp}{2kh}\right) \Rightarrow \cot (2\phi )=\dfrac{kq-hp}{2kh}%
\Rightarrow \dfrac{\cot ^{2}\phi -1}{2\cot \phi }=\dfrac{kq-hp}{2kh}%
\Rightarrow \cot \phi =\dfrac{1}{2kh}\left( kq-hp\pm \sqrt{%
4k^{2}h^{2}+(kq-hp)^{2}}\right) =\dfrac{kq-hp\pm \sqrt{W}}{2kh}$. We first
need to determine whether to choose the positive or the negative root. If $%
kq-hp>0$, then $\cot (2\phi )=\dfrac{kq-hp}{2kh}>0\Rightarrow 0\leq 2\phi
\leq \dfrac{\pi }{2}\Rightarrow 0\leq \phi \leq \dfrac{\pi }{4}\Rightarrow
1<\cot \phi <\infty $. Let $x=2kh$, $y=kq-hp$, $0<x<\infty $, $0<y<\infty $.
If $\cot \phi =\dfrac{kq-hp-\sqrt{W}}{2kh}$, then $\cot \phi =\dfrac{y-\sqrt{%
x^{2}+y^{2}}}{x}=\dfrac{y}{x}-\sqrt{1+\left( \dfrac{y}{x}\right) ^{2}}=u-%
\sqrt{1+u^{2}},$ where $u=\dfrac{y}{x}$, $0<u<\infty $. Let $z(u)=u-\sqrt{%
1+u^{2}}$. Then $z^{\prime }(u)=\allowbreak \dfrac{\sqrt{1+u^{2}}-u}{\sqrt{%
1+u^{2}}}>0$, $z(0)=\allowbreak -1$, and $\lim\limits_{u\rightarrow \infty
}z(u)=0$. Thus $-1<z(u)<0\Rightarrow -1<\cot \phi <0$, which contradicts $%
1<\cot \phi <\infty $. If $kq-hp<0,$ then $\cot (2\phi )=\dfrac{kq-hp}{2kh}%
<0\Rightarrow \dfrac{\pi }{2}\leq 2\phi \leq \pi \Rightarrow \dfrac{\pi }{4}%
\leq \phi <\dfrac{\pi }{2}\Rightarrow 0<\cot \phi <1$. Again, if $\cot \phi =%
\dfrac{kq-hp-\sqrt{W}}{2kh}$, then $\cot \phi =z(u),-\infty <u<0$. Since $%
z(0)=\allowbreak -1$ and $\lim\limits_{u\rightarrow -\infty
}z(u)=\allowbreak -\infty $, $-\infty <z(u)<-1\Rightarrow \cot \phi <-1$,
which contradicts $0<\cot \phi <1$. That proves 
\begin{equation}
\cot \phi =\dfrac{kq-hp+\sqrt{W}}{2kh}.  \tag{cotphi}
\end{equation}

To finish the proof of Proposition 2, note that $M_{1}=\dfrac{-kh+\sqrt{kh}%
\sqrt{kh-pq}}{hp}=\sqrt{kh}\dfrac{-\sqrt{kh}+\sqrt{kh-pq}}{hp}<0$ and $%
M_{2}=-\dfrac{k}{p}-\sqrt{kh}\dfrac{\sqrt{kh-pq}}{hp}<0$. Thus the only way
that $L$ and $L^{\prime }$ can form angles of $\theta $ and $-\theta $,
respectively, with the major axis of $E_{O}$ is if the major axis of $E_{O}$
has a negative slope. In that case the \textbf{minor} axis of $E_{O}$ is
rotated by $\phi $ counterclockwise from the positive $x$ axis. It follows
that the two directional constants, $M$ and $M^{\prime }$, are given by $%
\tan \left( \phi +\theta -\dfrac{\pi }{2}\right) $ and $\tan \left( \phi
-\theta -\dfrac{\pi }{2}\right) $. We shall prove that $\tan \left( \phi
+\theta -\dfrac{\pi }{2}\right) =M_{1}$. We find it convenient to introduce
the following variables: 
\begin{equation*}
s=\allowbreak hp+kq,t=\allowbreak hp-kq.
\end{equation*}

Note that $2k^{2}h-kpq+hp^{2}=k(kh-pq)+k^{2}h+hp^{2}>0$ by (1). Hence $%
(kv_{0}+h)+\sqrt{(kv_{0}-h)^{2}+\left( v_{0}p+q\right) ^{2}}=kv_{0}+h+\dfrac{%
\left( ph+qk\right) \sqrt{W}}{2k^{2}h-kpq+hp^{2}}$, which implies that

$\dfrac{(kv_{0}+h)\left( 2k^{2}h-kpq+hp^{2}\right) +\left( ph+qk\right) 
\sqrt{W}}{2k^{2}h-kpq+hp^{2}}=\dfrac{W+\left( ph+qk\right) \sqrt{W}}{%
2k^{2}h-kpq+hp^{2}}=$

$\sqrt{W}\dfrac{\sqrt{W}+\left( ph+qk\right) }{2k^{2}h-kpq+hp^{2}}$. By (fv)
and (gv0), $f(v_{0})=\dfrac{4kh\left( kh-pq\right) W}{\left(
2k^{2}h-kpq+hp^{2}\right) ^{2}}\dfrac{\left( 2k^{2}h-kpq+hp^{2}\right) ^{2}}{%
W\left( \sqrt{W}+\left( ph+qk\right) \right) ^{2}}=\dfrac{4kh\left(
kh-pq\right) }{\left( \sqrt{W}+\left( ph+qk\right) \right) ^{2}}=\dfrac{4r}{%
\left( \sqrt{W}+s\right) ^{2}}$. By (fv) again, 
\begin{equation}
\dfrac{b}{a}=\dfrac{2r}{\sqrt{W}+s}.  \tag{2}
\end{equation}

By (cotphi) and (2), $\tan \left( \phi +\theta -\dfrac{\pi }{2}\right) =%
\dfrac{\tan \theta \tan \phi -1}{\tan \theta +\tan \phi }=\dfrac{\dfrac{b}{a}%
\dfrac{2kh}{kq-hp+\sqrt{W}}-1}{\dfrac{b}{a}+\dfrac{2kh}{kq-hp+\sqrt{W}}}=%
\dfrac{\dfrac{2r}{\sqrt{W}+s}\dfrac{2kh}{\sqrt{W}-t}-1}{\dfrac{2r}{\sqrt{W}+s%
}+\dfrac{2kh}{\sqrt{W}-t}}=$

$\dfrac{4khr-\left( \sqrt{W}+s\right) \left( \sqrt{W}-t\right) }{2r\left( 
\sqrt{W}-t\right) +2kh\left( \sqrt{W}+s\right) }$. Hence $\tan \left( \phi
+\theta -\dfrac{\pi }{2}\right) -M_{1}=\dfrac{1}{2}\dfrac{4khr-\left( \sqrt{W%
}+s\right) \left( \sqrt{W}-t\right) }{r\left( \sqrt{W}-t\right) +kh\left( 
\sqrt{W}+s\right) }-\dfrac{r-hk}{hp}=0\iff $

$4kh^{2}rp-hp\left( \sqrt{W}+s\right) \left( \sqrt{W}-t\right)
-2r(r-hk)\left( \sqrt{W}-t\right) -2(r-hk)kh\left( \sqrt{W}+s\right) =0$

$\iff 4kh^{2}rp+hpst+2r(r-hk)t-2s(r-hk)kh+\left( -hp\left( s-t\right)
-2r(r-hk)-2(r-hk)kh\right) \sqrt{W}$

$-hpW=0$. Now $4kh^{2}rp+hpst+2r(r-hk)t-2s(r-hk)kh=hpW$ and

$\left( -hp\left( s-t\right) -2r(r-hk)-2(r-hk)kh\right) =0$. Hence $\tan
\left( \phi +\theta -\dfrac{\pi }{2}\right) =M_{1}$. Similarly, one can show
that $\tan \left( \phi -\theta -\dfrac{\pi }{2}\right) =M_{2}$. 
\endproof%

By Propositions 1 and 2 and the main result in [D], we have

\textbf{Theorem 1:} There exists a unique ellipse, $E_{O}$, which passes
through the vertices of $R$\ and whose \textit{equal} conjugate diameters
possess the directional constants $M_{1}$ and $M_{2}$. Furthermore, $E_{O}$
is the unique ellipse of minimal eccentricity among all ellipses which pass
through the vertices of $R$.

\section{Minimal Area}

We now prove a result similar to Proposition 1, but which instead minimizes
the \textbf{area }among all ellipses which pass through the vertices of $R$.
This was not discussed by Steiner in [D] and there does not appear to be a
nice characterization of the minimal area ellipse. Again we shall prove the
case when $R$\ is not a trapezoid. Since convex quadrilaterals and ratios of
areas of ellipses are preserved under one--one affine transformations, we
may assume throughout this section that the vertices of $R$\ are $%
(0,0),(1,0),(0,1)$, and $(s,t)$ for some positive real numbers $s$ and $t$.
Furthermore, since $R$\ is convex and not a trapezoid, it follows easily
that 
\begin{equation}
s+t>1\text{ and }s\neq 1\neq t.  \tag{3}
\end{equation}%
\textbf{Lemma 2:} Suppose that the vertices of $R$\ are $(0,0),(1,0),(0,1)$,
and $(s,t)$ for some positive real numbers $s$ and $t$ satisfying (3). Let $%
m_{s,t}=\dfrac{t}{s(s-1)^{2}}\left( s+t-1+st-2\sqrt{st\left( s+t-1\right) }%
\right) $ and $M_{s,t}=\dfrac{t}{s(s-1)^{2}}\left( s+t-1+st+2\sqrt{st\left(
s+t-1\right) }\right) $. An ellipse, $E_{0}$, passes through the vertices of 
$R$\ if and only if $E_{0}$ has the form

\begin{equation}
stux^{2}+sty^{2}-\left( s(s-1)u+t(t-1)\right) xy-stux-sty=0,u\in
I_{s,t}=\left( m_{s,t},M_{s,t}\right) .  \tag{5}
\end{equation}

If $a$ and $b$ denote the lengths of the semi--major and semi--minor axes,
respectively, of $E_{0}$, then%
\begin{equation}
a^{2}=\tfrac{2u(su+t)s^{2}t^{2}\left( s+t-1\right) }{\left( 4st\left(
s+t-1\right) u-\left( s(s-1)u-t(t-1)\right) ^{2}\right) \left( st(u+1)-\sqrt{%
t^{2}\left( s^{2}+(t-1)^{2}\right) -2st\left( s+t-1\right) \allowbreak
u+s^{2}\left( t^{2}+(s-1)^{2}\right) u^{2}}\right) }  \tag{6}
\end{equation}%
and 
\begin{equation}
b^{2}=\tfrac{2u(su+t)s^{2}t^{2}\left( s+t-1\right) }{\left( 4st\left(
s+t-1\right) u-\left( s(s-1)u-t(t-1)\right) ^{2}\right) \left( st(u+1)+\sqrt{%
t^{2}\left( s^{2}+(t-1)^{2}\right) -2st\left( s+t-1\right) \allowbreak
u+s^{2}\left( t^{2}+(s-1)^{2}\right) u^{2}}\right) }.  \tag{7}
\end{equation}

Finally, the center of $E_{0},\left( x_{0},y_{0}\right) ,$ is given by 
\begin{equation}
x_{0}=\dfrac{st\left( \left( 2st+s^{2}-s\right) u+(t^{2}-t)\right) }{%
2st\left( st+s+t-1\right) u-\allowbreak s^{2}\left( s-1\right)
^{2}u^{2}-t^{2}\allowbreak \left( t-1\right) ^{2}}  \tag{8}
\end{equation}%
and 
\begin{equation}
y_{0}=\dfrac{st\left( \left( 2st+t^{2}-t\right) u+(s^{2}-s)u^{2}\right) }{%
2st\left( st+s+t-1\right) u-\allowbreak s^{2}\left( s-1\right)
^{2}u^{2}-t^{2}\allowbreak \left( t-1\right) ^{2}}  \tag{9}
\end{equation}

\textbf{Proof: }Substituting the vertices of $R$\ into the general equation
of a conic, $Ax^{2}+By^{2}+2Cxy+Dx+Ey+F=0,$ yields the equations $F=0$, $%
A+D=0,B+E=0$, and $As^{2}+Bt^{2}+2Cst-As-Bt=0$, which implies that $%
As(s-1)+Bt(t-1)+2Cst=0$ or $C=-\dfrac{As(s-1)+Bt(t-1)}{2st}$. Multiplying
thru by $st$ and dividing thru by $B$ yields the equation in (5), with $u=%
\dfrac{A}{B}$. Conversely, any curve satisfying (5) must pass through the
vertices of $R$. The curve defined by $Ax^{2}+By^{2}+2Cxy+Dx+Ey+F=0$ is an
ellipse if and only if $AE^{2}+BD^{2}+4FC^{2}-2CDE-4ABF\neq 0$ and $%
AB-C^{2}>0$. The first condition becomes $s^{2}t^{2}u\left( s+t-1\right)
\left( su+t\right) \neq 0$, which clearly holds since $s+t>1$. The second
condition becomes

$4\left( stu\right) \left( st\right) -\left( s(s-1)u+t(t-1)\right)
^{2}=-s^{2}\left( s-1\right) ^{2}u^{2}+2st\left( st+s+t-1\right)
u-t^{2}\left( t-1\right) ^{2}>0$, which we write as $\alpha (u)<0$, where

\begin{equation*}
\alpha (u)=s^{2}(s-1)^{2}u^{2}-2st\left( st+s+t-1\right) u+t^{2}(t-1)^{2}.
\end{equation*}

Now it is easy to show that $\alpha (u)<0\iff m_{s,t}<u<M_{s,t}$. That
proves (5). If $E_{0}$ satisfies (5), then (6) and (7) follow immediately
from Lemma 1--(asq) and (bsq), and (8) and (9) follow immediately from Lemma
1--center, with $A=stu,B=st,C=-\dfrac{1}{2}\left( s(s-1)u+t(t-1)\right)
,D=-stu,E=-st,$ and $F=0$.

$\allowbreak $\textbf{Theorem 2:} There exists a unique ellipse, $E_{A}$, of
minimal area which passes through the vertices of $R$.

\textbf{Proof: }By Lemma 2--(6) and (7), $a^{2}b^{2}=\tfrac{%
4u^{2}(su+t)^{2}s^{2}t^{2}\left( st\left( s+t-1\right) \right) ^{2}}{\left(
4st\left( s+t-1\right) u-\left( s(s-1)u-t(t-1)\right) ^{2}\right) ^{2}}%
\times $

$\tfrac{1}{\left( st(u+1)-\sqrt{t^{2}\left( s^{2}+(t-1)^{2}\right)
-2st\left( s+t-1\right) \allowbreak u+\allowbreak s^{2}\left(
t^{2}+(s-1)^{2}\right) u^{2}}\right) }\times $

$\tfrac{1}{\left( st(u+1)+\sqrt{t^{2}\left( s^{2}+(t-1)^{2}\right)
-2st\left( s+t-1\right) u+s^{2}\left( t^{2}+(s-1)^{2}\right) u^{2}}\right) }=
$

$=\tfrac{4u^{2}(su+t)^{2}s^{2}t^{2}\left( st\left( s+t-1\right) \right) ^{2}%
}{\left( 4st\left( s+t-1\right) u-\left( s(s-1)u-t(t-1)\right) ^{2}\right)
^{2}\left( s^{2}t^{2}(u+1)^{2}-\left( t^{2}\left( s^{2}+(t-1)^{2}\right)
-2st\left( s+t-1\right) u+s^{2}\left( t^{2}+(s-1)^{2}\right) u^{2}\right)
\right) }$

$=\dfrac{4u^{2}(su+t)^{2}s^{2}t^{2}\left( st\left( s+t-1\right) \right) ^{2}%
}{\left( -t^{2}\left( t-1\right) ^{2}+\left( 4s^{2}t^{2}-2s\left( s-1\right)
t\left( t-1\right) \right) u-s^{2}\left( s-1\right) ^{2}u^{2}\right) ^{3}}%
=\beta (u)$, where

\begin{equation*}
\beta (u)=-\dfrac{4u^{2}(su+t)^{2}s^{2}t^{2}\left( st\left( s+t-1\right)
\right) ^{2}}{\left( \alpha (u)\right) ^{3}}.
\end{equation*}

Note that $\beta $ is differentiable on $I_{s,t}$ since $\alpha (u)<0$
there. Also, $m_{s,t}>0\iff s+t-1+st>2\sqrt{st\left( s+t-1\right) }\iff
\left( s+t-1+st\right) ^{2}>4st\left( s+t-1\right) $(since $s+t>1$) $\iff
\left( t-1\right) ^{2}\left( s-1\right) ^{2}>0$, which holds since $s,t\neq
1 $. Thus $m_{s,t}>0$ and $M_{s,t}>0$, which implies that $I_{s,t}\subset
\left( 0,\infty \right) $. Now $\lim\limits_{u\rightarrow m_{s,t}^{+}}\alpha
(u)=\lim\limits_{u\rightarrow M_{s,t}^{-}}\alpha (u)=0$ thru negative
numbers and the numerator of $\beta (u),$ for given $s$ and $t$, satisfies $%
4u^{2}(su+t)^{2}s^{2}t^{2}\left( st\left( s+t-1\right) \right) ^{2}>\delta
>0 $ for $u\in I_{s,t}$. Thus $\lim\limits_{u\rightarrow m_{s,t}^{+}}\beta
(u)=\lim\limits_{u\rightarrow M_{s,t}^{-}}\beta (u)=\infty $, which implies
that $\beta $ must attain its global minimum on $I_{s,t}$. $\beta ^{\prime
}(u)=\allowbreak 8u\left( su+t\right) s^{2}t^{2}\left( st\left( s+t-1\right)
\right) ^{2}\dfrac{\gamma (u)}{\left( \alpha (u)\right) ^{4}}$, where

\begin{equation*}
\gamma (u)=s^{3}\left( s-1\right) ^{2}u^{3}+s^{2}t\left(
2s^{2}-3s+st+1+t\right) \allowbreak u^{2}\allowbreak -st^{2}\left(
2t^{2}+st-3t+s+1\right) u-t^{3}\left( t-1\right) ^{2}.
\end{equation*}

$\gamma (0)=-t^{3}\left( t-1\right) ^{2}<0$ and $\lim\limits_{u\rightarrow
\infty }\gamma (u)=\infty $, which implies that $\gamma $ has at least one
real root in $\left( 0,\infty \right) $. We now show that $\gamma $ is
convex on $I_{s,t}$ by looking at $\gamma ^{\prime \prime }(u)=2s^{2}d(s,t)$%
, where $d(s,t)=2s^{2}+st-3s+t+1$. Now $\dfrac{\partial d}{\partial s}=4s+t-3
$ and $\dfrac{\partial d}{\partial t}=s+1$, which implies that $d$ has no
critical points in $S=\left\{ (s,t):s+t\geq 1,s\geq 0,t\geq 0\right\} $. The
boundary of $S$ consists of $S_{1}\cup S_{2}\cup S_{3},$ where $%
S_{1}=\left\{ (s,t):s+t=1,0\leq s\leq 1\right\} ,S_{2}=\left\{
(s,t):s=0,t\geq 1\right\} ,$ and $S_{3}=\left\{ (s,t):t=0,s\geq 1\right\} $.
We now check $d$ on $\partial \left( S\right) $. On $S_{1}$ we have $%
d(1-t,t)=t\left( 1+t\right) \geq 0$, on $S_{2}$ we have $d(0,t)=1+t\geq 0$,
and on $S_{3}$ we have $d(s,0)=\left( 2s-1\right) \left( s-1\right) \geq 0$.
Thus $d(s,t)\geq 0$ on $S$, which implies that $\gamma $ is convex on $%
I_{s,t}$. That in turn implies that $\beta $ has a \textbf{unique} global
minimum on $I_{s,t}$, which yields a unique ellipse of minimal area which
passes through the vertices of $R$.

\textbf{Remark: }In [RP] \& [RP2], the authors investigate the problem of
constructing and characterizing an ellipse of minimal area containing a
finite set of points. The results and methods in this paper are different
than those in [RP] \& [RP2], but it is worth pointing out some of the small
intersection. In particular, for a convex quadrilateral, $R$, the authors in
[RP] \& [RP2] construct an algorithm for finding the minimal area ellipse
containing $R$\ and they also prove a uniqueness result. For the case when
this ellipse passes thru all four vertices of $R$, this ellipse is then the
minimal area ellipse discussed in this paper. However, there are convex
quadrilaterals, $R$, for which the minimal area ellipse\ containing $R$\
does not pass thru all four vertices of $R$. In that case, the the minimal
area ellipse discussed in this paper is not the same.

\section{Inscribed versus Circumscribed}

In this section and the next, we allow $R$\ to be a \textbf{trapezoid}, so
we shall need a version of Lemma 2 for trapezoids. The proof of Lemma 3
follows immediately from Lemma 1 or from Lemma 2 by allowing $s$ to approach 
$1$. We omit the details here.

\textbf{Lemma 3:} Suppose that $R$\ is a \textbf{trapezoid }with vertices $%
(0,0),(1,0),(0,1)$, and $(1,t),t\neq 1$.

An ellipse, $E_{0}$, passes through the vertices of $R$\ if and only if $%
E_{0}$ has the form

\begin{equation}
tux^{2}+ty^{2}-t(t-1)xy-tux-ty=0,u\in I_{t}=\left( \dfrac{1}{4}\left(
t-1\right) ^{2},\infty \right)   \tag{10}
\end{equation}

If $a$ and $b$ denote the lengths of the semi--major and semi--minor axes,
respectively, of $E_{0}$, then 
\begin{equation}
a^{2}=\allowbreak \dfrac{-2u\left( u+t\right) }{\left( -4u+\left( t-1\right)
^{2}\right) \left( u+1-\sqrt{2+t^{2}-2t-2u+u^{2}}\right) }  \tag{11}
\end{equation}%
$\allowbreak $

and 
\begin{equation}
b^{2}=\allowbreak \dfrac{-2u\left( u+t\right) }{\left( \left( t-1\right)
^{2}-4u\right) \left( u+1+\sqrt{\left( t-1\right) ^{2}+\left( u-1\right) ^{2}%
}\right) }  \tag{12}
\end{equation}

Finally, the center of $E_{0},\left( x_{0},y_{0}\right) ,$ is given by 
\begin{equation}
x_{0}=\dfrac{2u+t-1}{4u-(t-1)^{2}},y_{0}=\dfrac{\left( 1+t\right) u}{%
4u-(t-1)^{2}}  \tag{13}
\end{equation}

\textbf{Theorem 3: }Let $R$\ be a convex quadrilateral in the $xy$ plane
which is \textbf{not} a parallelogram. Suppose that $E_{1}$ and $E_{2}$ are
each ellipses, with $E_{1}$ inscribed in $R$\ and $E_{2}$ circumscribed
about $R$. Then $E_{1}$ and $E_{2}$ cannot have the same center.

\textbf{Proof: }Assume first that $R$\ is \textbf{not} a \textbf{trapezoid}.
Since the center of an ellipse is affine invariant, we may assume that the
vertices of $R$\ are $(0,0),(1,0),(0,1)$, and $(s,t)$ as above, where $s$
and $t$ satisfy (3). By [H, Theorem 2.3], if $M_{1}$ and $M_{2}$ are the
midpoints of the diagonals of $R$, then each point on the open line segment, 
$Z$, connecting $M_{1}$ and $M_{2}$ is the center of some ellipse inscribed
in $R$. Thus the locus of centers of $E_{1}$ is precisely $Z$. For $R$\
above, the equation of $Z$ is $y=L(x)=\dfrac{1}{2}\dfrac{s-t+2x(t-1)}{s-1},$
where $x$ lies in the open interval connnecting $\dfrac{1}{2}$ and $\dfrac{1%
}{2}s$. If $E_{1}$ and $E_{2}$ have the same center, then the center of $%
E_{2},\left( x_{0},y_{0}\right) ,$ must lie on $Z$. Hence $L\left(
x_{0}\right) =y_{0}$, which implies that $L\left( x_{0}\right) -y_{0}=\dfrac{%
-\left( t+s\right) \left( \left( s-s^{2}\right) u+t^{2}-t\right) \left(
\left( s^{2}-s\right) u+t^{2}-t\right) }{2\left(
-2us^{2}t^{2}-2us^{2}t-2ust^{2}+2ust+s^{4}u^{2}-2s^{3}u^{2}+s^{2}u^{2}+t^{4}-2t^{3}+t^{2}\right) \left( s-1\right) 
}=0$. Thus $\allowbreak \left( s-s^{2}\right) u+t^{2}-t=0$ or $\left(
s^{2}-s\right) u+t^{2}-t=0$, which implies that $u=\pm \dfrac{t^{2}-t}{%
s^{2}-s}$. If $u=\dfrac{t^{2}-t}{s^{2}-s}$, then some simplification yields,
by (8) in Lemma 2, $x_{0}=\allowbreak \dfrac{1}{2}s$. Similarly, if $u=-%
\dfrac{t^{2}-t}{s^{2}-s}$, then $x_{0}=\dfrac{1}{2}$. But $\dfrac{1}{2}s$
and $\dfrac{1}{2}$ do not lie on $Z$, and thus $E_{1}$ and $E_{2}$ cannot
have the same center. Now suppose that $R$\ is a trapezoid, but not a
parallelogram. Then we may assume, again by affine invariance, that the
vertices of $R$\ are $(0,0),(1,0),(0,1)$, and $(1,t),t\neq 1$. The equation
of $Z$ is now $x=\dfrac{1}{2},$ where $y$ lies in the open interval
connnecting $\dfrac{1}{2}$ and $\dfrac{1}{2}t$. If $E_{1}$ and $E_{2}$ have
the same center, then $x_{0}=\dfrac{1}{2}$. By (13) of Lemma 3, $\dfrac{%
2u+t-1}{4u-(t-1)^{2}}=\dfrac{1}{2}\Rightarrow
4u+2t-2=4u-(t-1)^{2}\Rightarrow t=\pm 1$, which contradicts the assumption
that $t>0,t\neq 1$. Hence again $E_{1}$ and $E_{2}$ cannot have the same
center.

It is easy to find examples where the center of an ellipse circumscribed
about $R$\ may lie inside, on, or outside $R$. We make the following
conjectures.

\textbf{Conjecture 1:} The center of the ellipse of minimal eccentricity
circumscribed about $R$\ lies inside $R$.

\textbf{Conjecture 2:} The center of the ellipse of minimal area
circumscribed about $R$\ lies inside $R$.

\section{Bielliptic Quadrilaterals}

The following definition is well--known.

\textbf{Definition 1: }Let\textbf{\ }$R$\ be a convex quadrilateral in the $%
xy$ plane.

(A) $R$\ is called cyclic if there is a circle which passes through the
vertices of $R$.

(B) $R$\ is called tangential if a circle can be inscribed in $R$.

(C) $R$\ is called bicentric if $R$\ is both cyclic and tangential.

\qquad We generalize the notion of bicentric quadrilaterals as follows. In
[H, Theorem 4.4] the author proved that there is a unique ellipse, $E_{I}$,
of minimal eccentricity inscribed in a convex quadrilateral, $R$. Using
Proposition 1 from this paper, we let $E_{O}$ be the unique ellipse of
minimal eccentricity circumscribed about $R$.

\textbf{Definition 2:} A convex quadrilateral is called \textbf{bielliptic}
if $E_{I}$ and $E_{O}$ have the \textbf{same eccentricity}.

If $R$\ is bielliptic, we say that $R$\ is of class $\tau $, $0\leq \tau <1$%
, if $E_{I}$ and $E_{O}$ each have eccentricity $\tau $.

It is natural to ask the following:

\textbf{Question:} Does there exist a bielliptic quadrilateral of class $%
\tau $ for \textit{some} $\tau ,\tau >0$ ?

We answer this in the affirmative with the following results.

\textbf{Theorem 4: }There exists a convex quadrilateral, $R$, which is
bielliptic of class $\tau $ for \textit{some} $\tau >0$. That is, there
exists a bielliptic convex quadrilateral which is not bicentric.

\textbf{Proof: }Consider the convex quadrilateral, $R$, with vertices $%
(0,0),(1,0),(0,1)$, and $(s,t)$. We shall show that for some $s$ and $t$
satisfying (3), $R$\ is bielliptic of class $\tau $ for some $\tau >0$. It
is easy to show that $R$\ is cyclic if and only if $\left( s-\dfrac{1}{2}%
\right) ^{2}+\left( t-\dfrac{1}{2}\right) ^{2}=\dfrac{1}{2}$. In general, a
convex quadrilateral is tangential if and only if opposite sides add up to
the same sum. It follows that $R$\ is tangential if and only if $s=t$. Thus $%
s=\dfrac{1}{2}$, $t=\dfrac{1}{2}\left( 1+\sqrt{2}\right) \approx \allowbreak
1.\allowbreak 207$ gives a cyclic quadrilateral, and $s=2$, $t=2$ gives a
tangential quadrilateral. Consider the family of quadrilaterals $R_{r}$
given by 
\begin{equation}
s=-\dfrac{3}{2}r+2,t=r\left( \dfrac{1}{2}+\dfrac{1}{2}\sqrt{2}\right)
+2-2r,0\leq r\leq 1.  \tag{14}
\end{equation}%
It is easy to see that $r=0$\textbf{\ }gives a tangential quadrilateral
which is not cyclic, and $r=1$\ gives a cyclic quadrilateral which is not
tangential. Since the eccentricity of the inscribed and circumscribed
ellipses of minimal eccentricity, $E_{I}(r)$ and $E_{O}(r)$, each vary
continuously with $r$, $R_{r}$ must be bielliptic for some $r,0<r<1$. More
precisely, let $\epsilon _{I}(r)$ and $\epsilon _{O}(r)$ denote the
eccentricities of $E_{I}$ and $E_{O}$, respectively. Then $\epsilon _{I}(0)=0
$ and $\epsilon _{O}(0)>0$ since $E_{I}(0)$ is a circle and $E_{O}(0)$ is
not a circle. Similarly, $\epsilon _{I}(1)>0$ and $\epsilon _{O}(1)=0$ since 
$E_{I}(1)$ is not a circle and $E_{O}(1)$ is a circle. Since $\epsilon
_{I}(r)$ and $\epsilon _{O}(r)$ are each continuous functions of $r,$ by the
Intermediate Value Theorem, $\epsilon _{I}(r_{0})=\epsilon _{O}(r_{0})$ for
some $0<r_{0}<1$. Now if $s$ and $t$ satisfy (14), then $%
(2s-1)^{2}+(2t-1)^{2}-2=-\dfrac{2}{41}\left( -10+3\sqrt{2}\right) \left(
-1+r\right) \left( 41r-40-12\sqrt{2}\right) =0$

$\iff r=1$ or $r=\dfrac{40}{41}+\dfrac{12}{41}\sqrt{2}\approx \allowbreak
1.\allowbreak 39>1$. So for $0<r<1$, $R_{r}$ cannot be cyclic. Also, $%
s=t\iff -\dfrac{3}{2}r+2=-\dfrac{3}{2}r+\dfrac{1}{2}r\sqrt{2}+2\iff r=0$. So
for $0<r<1$, $R_{r}$ cannot be tangential. It follows that $\epsilon
_{I}(r_{0})=\epsilon _{O}(r_{0})=\tau >0$, which menas that $R_{r_{0}}$ is
bielliptic of class $\tau $.

\textbf{Theorem 5: }There exists a \textit{bielliptic trapezoid} of class $%
\tau $ for \textit{some} $\tau >0$.

\textbf{Proof: }Consider the trapezoid, $R$, with vertices $(0,0),(1,0),(0,1)
$, and $(1,t),t\neq 1$. We shall show that for some $t\neq 1,$ $R$\ is
bielliptic of class $\tau >0$. By Lemma 3--(11) and (12), $\dfrac{b^{2}}{%
a^{2}}=\dfrac{\left( \left( t-1\right) ^{2}-4u\right) \left( u+1-\sqrt{%
\left( t-1\right) ^{2}+\left( u-1\right) ^{2}}\right) }{\left( \left(
t-1\right) ^{2}-4u\right) \left( u+1+\sqrt{\left( t-1\right) ^{2}+\left(
u-1\right) ^{2}}\right) }$. Hence the square of the eccentricity of an
ellipse circumscribed in $R$\ is given by $\epsilon (u)=1-\dfrac{b^{2}}{a^{2}%
}=\dfrac{2\sqrt{\left( t-1\right) ^{2}+\left( u-1\right) ^{2}}}{u+1+\sqrt{%
\left( t-1\right) ^{2}+\left( u-1\right) ^{2}}},u\in I_{t}=\left( \dfrac{1}{4%
}\left( t-1\right) ^{2},\infty \right) $. $\epsilon ^{\prime }(u)=\dfrac{%
-2\left( 3+t^{2}-2t-2u\right) }{\left( u+1+\sqrt{t^{2}-2t+2+u^{2}-2u}\right)
^{2}\sqrt{t^{2}-2t+2+u^{2}-2u}}=0\iff u=\dfrac{1}{2}\left( t^{2}-2t+3\right) 
$. We shall show that this value of $u$ gives the minimal eccentricity.
First, $\epsilon \left( \dfrac{1}{2}\left( t^{2}-2t+3\right) \right) =\dfrac{%
2\sqrt{\left( t^{2}-2t+5\right) \left( t-1\right) ^{2}}}{t^{2}-2t+5+\sqrt{%
\left( t^{2}-2t+5\right) \left( t-1\right) ^{2}}}=$

$\dfrac{2\left\vert t-1\right\vert \sqrt{t^{2}-2t+5}}{t^{2}-2t+5+\left\vert
t-1\right\vert \sqrt{t^{2}-2t+5}}=\dfrac{2\left\vert t-1\right\vert }{\sqrt{%
(t-1)^{2}+4}+\left\vert t-1\right\vert }<\dfrac{2\left\vert t-1\right\vert }{%
\left\vert t-1\right\vert +\left\vert t-1\right\vert }=\allowbreak 1$. Also, 
$\lim\limits_{u\rightarrow \left( \left( t-1\right) ^{2}/4\right)
^{+}}\epsilon (u)=\allowbreak 1$ and $\lim\limits_{u\rightarrow \infty
}\epsilon (u)=\allowbreak 1$. Thus the square of the minimal eccentricity of
an ellipse circumscribed in $R$\ is given by

\begin{equation}
\epsilon (t)=\allowbreak \dfrac{2\left\vert t-1\right\vert }{\sqrt{%
(t-1)^{2}+4}+\left\vert t-1\right\vert }  \tag{15}
\end{equation}

In [H] the author derived formulas for the eccentricity of the unique
ellipse of minimal eccentricity inscribed in a convex quadrilateral, $R$.
Those formulas apply when $R$\ is \textbf{not} a \textbf{trapezoid}. The
methods used in [H] can easily be adapted to the case when $R$\ is a
trapezoid. The ellipse of minimal eccentricity inscribed in a trapezoid is
also unique, and one can derive the following formulas. Let $I_{t}$ denote
the open interval with $\dfrac{1}{2}$ and $\dfrac{1}{2}t$ as endpoints. Let 
\begin{equation*}
E(k)=\dfrac{(2k-1)\left( 2k-t\right) }{16\left( t-1\right) ^{2}k^{4}+\left(
8+8t^{2}+48t\right) k^{2}-32t\left( 1+t\right) k+17t^{2}-2t+1},
\end{equation*}

\begin{equation}
\epsilon (k)=\dfrac{2}{1+\sqrt{1-16t\left( 1-t\right) ^{2}E(k)}},  \tag{16}
\end{equation}

and

\begin{equation*}
c(k)=16k^{3}-12(t+1)k^{2}+4(2t-1)k+t+1.
\end{equation*}

Then $c(k)$ has a unique root, $k_{0}$, in $I_{k}$, and $\epsilon \left(
k_{0}\right) $ equals the square of the minimal eccentricity of an ellipse
inscribed in $R$. By (15) and (16), we want to show that there is a value of 
$t\neq 1$ and $k\in I_{t}$ such that $c(k)=0$ and $\dfrac{2\left\vert
t-1\right\vert }{\sqrt{(t-1)^{2}+4}+\left\vert t-1\right\vert }=\dfrac{1}{1+%
\sqrt{1-16t\left( 1-t\right) ^{2}E(k)}}$. This is equivalent, after some
algebraic simplification, to $4t\left( t-1\right) ^{4}E(k)+1=0$. Some more
algebraic simplification yields the equation

\begin{gather}
16\left( t-1\right) ^{2}k^{4}+\left(
16t^{5}-64t^{4}+96t^{3}-56t^{2}+64t+8\right) k^{2}-  \tag{17} \\
8t\left( 1+t\right) \left( t^{2}-4t+5\right) \left( t^{2}+1\right) k+  \notag
\\
4t^{6}-16t^{5}+24t^{4}-16t^{3}+21t^{2}-2t+1=0  \notag
\end{gather}

Thus we want a solution to the simultaneous equations (17) and $c(k)=0$,
with $t\neq 1$ and $k\in I_{t}$. Maple gives the following solutions: $t=1,k=%
\dfrac{1}{2},k=\pm \dfrac{1}{2}i,t=\dfrac{1}{2}i,\allowbreak $ and $t=\dfrac{%
2\rho _{2}^{3}-3\rho _{2}^{2}+1-2\rho _{2}}{3\rho _{2}^{2}-4\rho _{2}-1},k=%
\dfrac{1}{2}\rho _{2}\allowbreak $ where $\rho _{2}$ is a root of

\begin{gather*}
p(x)=32x^{11}-287x^{10}+1006x^{9}-1487x^{8}+160x^{7}+ \\
1762x^{6}-884x^{5}-822x^{4}+80x^{3}+333x^{2}+150x+21
\end{gather*}

$t=1$ or $t=\dfrac{1}{2}i$ do not satisfy $t$ real, $t\neq 1$. Since $%
p(1)=\allowbreak 64>0$ and $p(1.5)=\allowbreak -23.\allowbreak 07715<0,$ $p$
must have a root, $x_{0}$, between $1$ and $2$. Numerically $x_{0}\approx
1.\allowbreak 232267$. It appears that the real roots of $p$ are
approximately $-0.8295535,1.\allowbreak 232267,1.778672$, though we don't
need that here. Now $\rho _{2}=1.\allowbreak 232267\Rightarrow t=\dfrac{%
2\rho _{2}^{3}-3\rho _{2}^{2}+1-2\rho _{2}}{3\rho _{2}^{2}-4\rho _{2}-1}%
\approx \allowbreak 1.\allowbreak 658119$. Then $k=\dfrac{1}{2}\rho
_{2}=0.6161335\in I_{t}$. The corresponding common value of the eccentricity
is $\approx \allowbreak 0.69013$.

\textbf{Remark:} It is interersting to note here that the bielliptic
quadrilateral in Theorem 4 is not a trapezoid. The family of quadrilaterals $%
R_{r}$ given in the proof of Theorem 4 yields a trapezoid if and only if $s=1
$ or $t=1$. Now $s=1\iff -\dfrac{3}{2}r+2=1\iff r=\dfrac{2}{3}$ and $t=1\iff
r\left( \dfrac{1}{2}+\dfrac{1}{2}\sqrt{2}\right) +2-2r=1\iff r=\dfrac{2}{3-%
\sqrt{2}}>1$. Thus $R_{r}$ is a trapezoid $\iff r=\dfrac{2}{3}$. Now $r=%
\dfrac{2}{3}\Rightarrow t=\allowbreak 1+\dfrac{1}{3}\sqrt{2}$. By (15) in
the proof of Theorem 5, the square of the minimal eccentricity of an ellipse
circumscribed about $R_{2/3}$\ is $\dfrac{2\left\vert t-1\right\vert }{\sqrt{%
(t-1)^{2}+4}+\left\vert t-1\right\vert }=\dfrac{2}{\sqrt{19}+1}\approx 0.373$%
. Also, $I_{t}=\left( \dfrac{1}{2},\dfrac{1}{2}t\right) \approx \allowbreak
\left( 0.5,0.736\right) $ and $c(k)=16k^{3}+\left( -24-4\sqrt{2}\right)
k^{2}+\left( \dfrac{8}{3}\sqrt{2}+4\right) k+2+\dfrac{1}{3}\sqrt{2}%
\allowbreak =0$ has the root $k\approx 0.5918015$ in $I_{t}$. That yields $%
E(k)\approx \allowbreak -1.\allowbreak 430$. By (16) in the proof of Theorem
5, the square of the minimal eccentricity of an ellipse inscribed in $R_{2/3}
$\ is $\epsilon ^{2}(k)\approx \allowbreak 0.511$. Thus the bielliptic
convex quadrilateral from Theorem 4 is not a trapezoid.

Theorems 4 and 5 show the existence of a bielliptic quadrilateral of class $%
\tau $ for some $0<\tau <1$. We cannot yet answer the following:

\textbf{Question:} Does there exist a bielliptic quadrilateral of class $%
\tau $ for \textit{each} $\tau $, $0<\tau <1$ ?

It would also be nice to answer

\textbf{Question: }If $R$\ is a bielliptic quadrilateral, is there a nice
relationship between the ellipse of minimal eccentricity inscribed in $R$\
and the ellipse of minimal eccentricity passing thru the vertices of $R$\ ?
This would generalize the known relationship between the inscribed and
circumscribed circles of bicentric quadrilaterals.

\section{References}

[D] Heinrich D\"{o}rrie, \textquotedblleft The most nearly circular ellipse
circumscribing a quadrilateral\textquotedblright , 100 Great problems of
Mathematics, Their History and Solution, Dover Publications, Inc. 1965,
231--236.

[H] Alan Horwitz, \textquotedblleft Ellipses of maximal area and of minimal
eccentricity inscribed in a convex quadrilateral\textquotedblright ,
Australian Journal of Mathematical Analysis and Applications, 2(2005), 1-12.

[RP] B. V. Rublev and Yu. I. Petunin, \textquotedblleft Minimum--Area
Ellipse Containing a Finite Set of Points. 1\textquotedblright , Ukrainian
Mathematical Journal, Vol. 50, No. 7, 1998, 1115--1124.

[RP2] B. V. Rublev and Yu. I. Petunin, \textquotedblleft Minimum--Area
Ellipse Containing a Finite Set of Points. 2\textquotedblright , Ukrainian
Mathematical Journal, Vol. 50, No. 8, 1998, 1253--1261.

[S] George Salmon, A treatise on conic sections, 6th edition, Chelsea
Publishing Company, New York.

\end{document}